\documentclass[12pt]{article}
\usepackage{amsmath,amssymb,amsthm,fullpage}

\title{A Tripos Surd}

\author{Mark B. Villarino\\
Depto.\ de Matem\'atica, Universidad de Costa Rica,\\
2060 San Jos\'e, Costa Rica}

\date{\today}

\theoremstyle{plain}
\newtheorem{thm}{Theorem}[section]    
\newtheorem{lemma}[thm]{Lemma}        

\theoremstyle{definition}

\numberwithin{equation}{section}

\makeatletter
\def\section{\@startsection{section}{1}{\z@}{-3.5ex plus -1ex minus
			  -.2ex}{2.3ex plus .2ex}{\large\bf}}
\def\subsection{\@startsection{subsection}{2}{\z@}{-3.25ex plus -1ex
			  minus -.2ex}{1.5ex plus .2ex}{\normalsize\bf}}
\renewcommand{\@dotsep}{200} 
\makeatother

\renewcommand{\geq}{\geqslant}  
\renewcommand{\leq}{\leqslant}  

\begin{document}

\maketitle

\section{Introduction}

By the final third of the nineteenth century the Mathematical Tripos examination had come to demand forty-four+ hours of grueling presssure-packed unrelenting high-level problem solving spread out over eight days.  This titanic tournament dwarfs the olympiads of today although it is legitimately regarded as their almost legendary progenitor.  Many of the problems set in this ``mega-olympiad"  were degree-level research questions and hundreds of papers have been written because of them.

Every year some tripos questions treated approximations and our interest was caught by the following striking fourth-root surd approximation taken from the tripos of 1886: 
\begin{quote} ``If $M=N^4+x$, and $x$ is small compared with $N$, then a good approximation for $\sqrt[4]{M}$ is:
\begin{equation}
\label{TA}
\frac{51}{56}N+\frac{5}{56}\frac{M}{N^3}+\frac{27Nx}{14(7M+5N^4)}.
\end{equation}Show that when $N=10,\ x=1$, this approximation is accurate to 16 places of decimals."
\footnote{We have not seen the examination, itself.  This quoted version is taken from Hardy \cite{Hardy}, p. 431.  The problem is also quoted by Chrystal \cite{Ch} p. 220 in a slightly different way...He writes ``p" instead of ``M" and leaves out ``compared with p," and writes ``approximately" instead of ``a good approximation ffor $\sqrt[4]{M}$ is..."  It is a shame that there is no online source for the problem statements in the old tripos examinations listed by year.}\end{quote}
Making the numerical substitutions we find, $$\sqrt[4]{10001}-\frac{1920160001}{192011200}=-5.695655...\cdot 10^{-18}$$
which shows that the approximation is in excess and just misses being accurate to $17$ decimal places!

Two questions naturally arise:\begin{enumerate}
  \item In general, how large is the \emph{error}:  \begin{equation}
\label{Ex}
E(x):=\sqrt[4]{M}-\left\{\frac{51}{56}N+\frac{5}{56}\frac{M}{N^3}+\frac{27Nx}{14(7M+5N^4)}\right\}?
\end{equation}
  \item How did the author \emph{discover} the formula \eqref{TA}?
\end{enumerate}

Answering these questions will lead us to some interesting and subtle mathematics.

\section{The Accuracy of the Surd Approximation}

We will prove
\begin{thm}  If $\frac{7}{12}\frac{|x|}{N^4}<1$, the error, $E(x)$,  is given by
\begin{equation}
\label{Ex2}
E(x)=-\frac{77}{2048}\left\{  \frac{1}{(1+X)^{\frac{15}{4}}}- \frac{28}{33}\frac{1}{1+\frac{7}{12}\frac{x}{N^4}}\right\}\frac{x^4}{N^{15}}
\end{equation}where $X$ is between $0$ and $\dfrac{x}{N^4}$. Moreover, if $$\frac{-20}{77}<\frac{x}{N^4}<.05336...$$the tripos surd approximation \eqref{TA} \textbf{overestimates} the true value of $\sqrt[4]{M}.$  
 \end{thm}If we apply \eqref{Ex2} to the original statement, and note that $0<X<\frac{1}{10000}$ we see that
 $$|E|\leq \frac{77}{2048}\left\{  \frac{1}{(1+0)^{\frac{15}{4}}}- \frac{28}{33}\frac{1}{1+\frac{7}{12}\frac{1}{10000}}\right\}\frac{1}{10^{15}}=5.698475...\cdot10^{-18}
$$which coincides with the true error up to $10^{-20}$ inclusive.

The form of the error term \eqref{Ex2} shows why the surd approximation is more accurate than the Taylor polynomial alone (see \eqref{M})...namely by subtracting the term $ \frac{28}{33}\frac{1}{1+\frac{7}{12}\frac{x}{N^4}}$ from the error  $\frac{1}{(1+X)^{\frac{15}{4}}}$ in the Taylor polynomial,  the approximation \emph{restores} part of the true sum lost by truncation.

 Although the formula \eqref{Ex2} for the error is exact, the presence of the unknown quantity $X$ can make it inconvenient in applications. Moreover, in order to transform \eqref{Ex2} into an inequality bounding $E(x)$ one needs information on the size of $x$.  Unfortunately, the original problem statement only says ``...$x$ is small compared with $N$..." which is not quantitatively precise.  Moreover, $x$ can be positive or negative which complicates the analysis.  We will prove
 
 \begin{thm}
 If $M$ differs from $N^4$ by less than $p\%$ of either, then $\sqrt[4]{M}$ differs from the tripos surd approximation by less than  
 $$\frac{77}{2048}\left\{\left(1+\frac{p}{100}\right)^{\frac{15}{4}}- \frac{28}{33}\right\}\left(\frac{p}{100+p}\right)^4\cdot N$$
  if $x$ is negative, and by less than $$\frac{77}{2048}\left\{ 1- \frac{28}{33}\frac{1}{1+\frac{7}{12}\frac{p}{100}}\right\}\left(\frac{p}{100}\right)^{4}\cdot N$$if $x$ is positive.
 
 \end{thm}
 
 Thus if $p=1$ then $\sqrt[4]{M}$ differs from the tripos surd approximation by less than $\dfrac{\text{N}}{17000000000}$ if the difference is positive, and by less than $\dfrac{\text{N}}{14600000000}$ if it is negative.
 \\
 
 The proofs of the two theorems are based on the standard Maclaurin expansion of $\sqrt[4]{1+\frac{x}{N^4}}$ with the Lagrange form of the remainder which we state as a separate lemma: 
 
 \begin{lemma}
  If $1+\frac{x}{N^4}\geq 0$ there exists a number $X$ between $0$ and $\frac{x}{N^4}$ such that the following expansion is valid:
 \begin{equation}
\label{M}
 \sqrt[4]{1+\frac{x}{N^4}}=1+\frac{1}{4}\frac{x}{N^4}-\frac{3}{32}\frac{x^2}{N^8}+\frac{7}{128}\frac{x^3}{N^{12}}-\frac{77}{2048}\left( \frac{1}{1+X}   \right)^{\frac{15}{4}}\frac{x}{N^4}.
  \end{equation}
  \end{lemma}
  \qed
 
  \begin{proof}[Proof of Theorem 2.1]
   Let $S(x)$ be the surd approximation:
  
  \begin{equation}
\label{Sx}
S(x):=\frac{51}{56}N+\frac{5}{56}\frac{M}{N^3}+\frac{27Nx}{14(7M+5N^4)}.
\end{equation}
Then, using $M=N^4+x$ we obtain
\begin{align*}
\label{}
    S(x)=& N\left\{\frac{51}{56}+\frac{5}{56}\frac{M}{N^4}+\frac{27x}{14(7M+5N^4)}.\right\}  \\
    =&  N\left\{\frac{51}{56}+\frac{5}{56}\frac{N^4+x}{N^4}+\frac{27x}{14(7M+5N^4)}.\right\}  \\
    =&N\left\{1+\frac{5}{56}\frac{x}{N^4}+\frac{27x}{14(7[N^4+x]+5N^4)}.\right\}  \\
    =&N\left\{1+\frac{5}{56}\frac{x}{N^4}+\frac{9}{56}\frac{x}{N^4}\frac{1}{1+\frac{7}{12}\frac{x}{N^4}}\right\}  \\
    =&N\left\{1+\frac{1}{4}\frac{x}{N^4}-\frac{3}{32}\frac{x^2}{N^8}+\frac{7}{128}\frac{x^3}{N^{12}}-\frac{49}{1536}\frac{x^4}{N^{16}}\frac{1}{1+\frac{7}{12}\frac{x}{N^4}}\right\}  
    \end{align*}
    where we have used the standard properties of the geometric series and the assumption $\frac{7}{12}\frac{|x|}{N^4}<1$.  But \eqref{M}  and $M=N^4+x$ show us that 
    $$\sqrt[4]{M}\equiv \sqrt[4]{N^4+x}=N\left\{1+\frac{1}{4}\frac{x}{N^4}-\frac{3}{32}\frac{x^2}{N^8}+\frac{7}{128}\frac{x^3}{N^{12}}-\frac{77}{2048}\left( \frac{1}{1+X}   \right)^{\frac{15}{4}}\frac{x}{N^4}\right\}$$Subtracting we obtain
    $$E(x):=\sqrt[4]{M}-S(x)=-\frac{77}{2048}\left\{  \frac{1}{(1+X)^{\frac{15}{4}}}- \frac{28}{33}\frac{1}{1+\frac{7}{12}\frac{x}{N^4}}\right\}\frac{x^4}{N^{15}}$$for some $X$ between $0$ and $\frac{x}{N^4}$.  This completes the proof of the formula.
    
    The proof that $S(x)$ overestimates $\sqrt[4]{1+\frac{x}{N^4}}$ is more troublesome because of the uncertainty of the value of $X$.  Analytically, we have to prove that for certain positive and negative values of $x$ the following inequality is valid:
  $$ \frac{1}{(1+X)^{\frac{15}{4}}}- \frac{28}{33}\frac{1}{1+\frac{7}{12}\frac{x}{N^4}}>0$$In case $0<X<\frac{x}{N^4}$ then
  $$\frac{1}{(1+X)^{\frac{15}{4}}}- \frac{28}{33}\frac{1}{1+\frac{7}{12}\frac{x}{N^4}}>\frac{1}{(1+\frac{x}{N^4})^{\frac{15}{4}}}- \frac{28}{33}\frac{1}{1+\frac{7}{12}\frac{x}{N^4}}$$and this latter function is positive for $\frac{x}{N^4}<.05336...$.  If $0>X>\frac{x}{N^4}$ then
   $$\frac{1}{(1+X)^{\frac{15}{4}}}- \frac{28}{33}\frac{1}{1+\frac{7}{12}\frac{x}{N^4}}>1- \frac{28}{33}\frac{1}{1+\frac{7}{12}\frac{x}{N^4}}$$which is positive for $\frac{x}{N^4}>\frac{-20}{77}=-.259...$.  Therefore
   $$\frac{-20}{77}<\frac{x}{N^4}<.05336...\Rightarrow\frac{1}{(1+X)^{\frac{15}{4}}}- \frac{28}{33}\frac{1}{1+\frac{7}{12}\frac{x}{N^4}}>0.$$Doubtless these bounds can be improved but now we have a quantitative formulation of ``...$x$ is small compared with $N$".  We note that the example in the original statement has $\frac{x}{N^4}=\frac{1}{10000}=.0001<.05336...$ which fulfills our inequality with more than plenty to spare..

  \end{proof}
  
  \begin{proof}[Proof of Theorem 2.2]:
  \\
  
  Suppose that $x$ is positive and $0<X<\frac{x}{N^4}$.  Then, by assumption
\begin{align*}
\label{}
 0<X<\frac{x}{N^4}<\frac{p}{100}\Rightarrow &1<1+X<1+\frac{x}{N^4} <1+\frac{p}{100}  \\
 \Rightarrow &1+\frac{7}{12}\frac{x}{N^4}<1+\frac{7}{12}\frac{p}{100}\\
  \Rightarrow &-\frac{1}{1+\frac{7}{12}\frac{x}{N^4}}<-\frac{1}{1+\frac{7}{12}\frac{p}{100}}\\
  \Rightarrow &|E(x)|:=\left|-\frac{77}{2048}\left\{  \frac{1}{(1+X)^{\frac{15}{4}}}- \frac{28}{33}\frac{1}{1+\frac{7}{12}\frac{x}{N^4}}\right\}\right|\cdot N\\
    <&  \frac{77}{2048}\left\{ 1- \frac{28}{33}\frac{1}{1+\frac{7}{12}\frac{p}{100}}\right\}\left(\frac{p}{100}\right)^{4}\cdot N
\end{align*}since $\dfrac{1}{(1+X)^{\frac{15}{4}}}<1.$
\\

Suppose that $x$ is negative and that $\frac{x}{N^4}<X<0.$ Then $M=N^4-|x|$ and
\begin{align*}
\label{}
   |x|<\frac{p}{100}\cdot M\Rightarrow & |x|<\frac{p}{100}\cdot (N^4-|x|)  \\
   \Rightarrow & \frac{|x|}{N^4} <\frac{p}{100}\left( 1-  \frac{|x|}{N^4}\right)\\
  \Rightarrow& \frac{|x|}{N^4}<\frac{p}{100+p}\\
  \Rightarrow&-\frac{p}{100+p}<\frac{x}{N^4}<X<0\\
  \Rightarrow&1-\frac{p}{100+p}<1+\frac{x}{N^4}<1+X<1\\
    \Rightarrow    &\frac{1}{1+X}<1+\frac{p}{100}\\
    \Rightarrow    &\frac{1}{(1+X)^{\frac{15}{4}}}<\left(1+\frac{p}{100}\right)^{\frac{15}{4}}.
   \end{align*}
   
   Moreover, 
   $$-\frac{1}{1+\frac{7}{12}\frac{x}{N^4}}<-1.$$Therefore, the formula \eqref{Ex2} the third and the last two inequalities above allow us to conclude that

 $$ |E(x)|<\frac{77}{2048}\left\{\left(1+\frac{p}{100}\right)^{\frac{15}{4}}- \frac{28}{33}\right\}\left(\frac{p}{100+p}\right)^4\cdot N$$This completes the proof.

  \end{proof}
  
  \section{Discovering the Approximation}
  
  We seek an approximation, $s(x)$, of the form
  
  \begin{equation}
\label{Sx2}
\sqrt[4]{M}\equiv\sqrt[4]{N^4+x}\approx AN+B\frac{M}{N^3}+\frac{CNx}{DM+EN^4}
\end{equation}where the coefficients $A,B,C,D,E$ are to be determined so that the approximation is as accurate as possible.  This means that it coincides with the Maclaurin expansion to as high a power as possible.

Expanding the right-hand side of $s(x)$, i.e., of \eqref{Sx2} into powers of $\frac{x}{N^4}$ we obtain,
\begin{equation}
\label{ }
s(x)=N\left\{A+B +\left(B+\frac{C}{D+E}\right)\frac{x}{N^4}-\frac{CD}{(D+E)^2}\frac{x^2}{N^8}+\frac{CD^2}{(D+E)^3}\frac{x^3}{N^{12}}-\cdots\right\}
\end{equation}Comparing this with \eqref{M} we obtain the following system of equations:
\begin{align*}
\label{}
  A+B=  &  1 \\
   B+\frac{C}{D+E}=&\frac{1}{4}  \\
   -\frac{CD}{(D+E)^2}=&-\frac{3}{32}\\
   \frac{CD^2}{(D+E)^3}=&\frac{7}{128}
\end{align*}The last two equations give us $E=\frac{5}{7}D$.  Then we obtain, in order, $\frac{C}{D+E}=\frac{5}{56}$, $B=\frac{5}{56}$,   $A=\frac{51}{56}$, $C=\frac{27}{98}D$.  Substituting the values of $E$ and $C$ into the fraction in $s(x)$ the common factor $D$ in the numerator and denominator cancel and it collapses to the fraction in $S(x)$.  This shows us that $s(x)\equiv S(x)$ and that $s(x)$ is uniquely determined.

The fact that four equations determine five unknowns is true and, indeed the values of $E$ and $C$ turn out to be multiples of a fifth unknown, D.  But, as we noted above, it cancels in the fraction after the substitution is made.

It is interesting to note that if we equate the coefficient of $\frac{x^4}{N^{15}}$ in the maclaurin expansion, namely $\frac{-77}{2048}$ with the corresponding coefficient in $s(x)$, namely $\frac{CD^3}{(D+E)^4}$ so as to obtain a fifth equation for the five unknowns, the third, fourth, and new fifth equations give the \emph{inconsistent} result $\frac{D}{D+E}=\frac{7}{12}=\frac{11}{16}.$  Therefore $S(x)$ is the best possible and unique approximation of the given form.
\subsubsection*{Acknowledgment}
Support from the Vicerrector\'{\i}a de Investigaci\'on of the 
University of Costa Rica is ack\-now\-ledged.

\end{document}